\documentstyle[11pt]{article}
\begin{document}

\def\b#1{{\bf #1}}
\def\i#1{{\it #1}}
\def\prob#1{{\bf P}\{#1\}}
\def\e#1{{\rm e}^{#1}}
\def\m#1{\mit #1}

\title{Sampling from a couple of positively correlated beta variates}
\author{Mario Catalani\\
Department of Economics\\ Via Po 53, 10124 Torino, Italy\\
E-mail: mario.catalani@unito.it}
\date{}
\maketitle
\begin{abstract}
\small{We know that the marginals in a Dirichlet distribution
are beta variates
exhibiting a negative correlation. But we can construct two linear combinations
of such marginals in such a way to obtain a positive correlation.
We discuss the restrictions that are to be imposed on the parameters
to accomplish such a result. In the case the sampling from the Dirichlet
distribution is performed through a generalization of Johnk's method
we discuss the efficiency of the algorithm implementing the method.
}
\end{abstract}

\section{Introduction}
Let $\{X_1,\,X_2,\,X_3\}$ be distributed according to a Dirichlet
distribution (see, for example, \cite[pages 231-235]{kotz}),
which will be denoted by $\{X_1,\,X_2,\,X_3\}\sim
{\cal D}(3,\,\alpha_0,\,\alpha_1,\,
\alpha_2,\,\alpha_3)$. To fix notation in this case
$$X_i={Z_i\over \sum_{j=0}^3 Z_j},\qquad i=1,\,2,\,3$$
where $Z_j,\;j=0,\,\ldots,\,3$ are independent \i{gamma} variables:
$Z_j\sim {\cal G}(\lambda,\,\alpha_j)$.

\noindent
Define
$$\left\{\begin{array}{lll}Y_1&=&X_1+X_3,\\
Y_2&=&X_2+X_3.\end{array}\right .$$
It follows
$$\left\{\begin{array}{lll}Y_1&\sim &{\cal B}e(\alpha_1+\alpha_3,\,
\alpha_0+\alpha_2)\\
Y_2&\sim &{\cal B}e(\alpha_2+\alpha_3,\,
\alpha_0+\alpha_1),
\end{array}\right .$$
where ${\cal B}e(\cdot,\,\cdot)$ denotes a \i{beta} distribution. Set $\gamma=
\alpha_0+\alpha_1+
\alpha_2+\alpha_3$. Then we have
\begin{eqnarray*}
&&\b{C}ov(X_1,\,X_2)=-{\alpha_1\alpha_2\over \gamma^2(\gamma+1)}\\
&&\b{C}ov(X_1,\,X_3)=-{\alpha_1\alpha_3\over \gamma^2(\gamma+1)}\\
&&\b{C}ov(X_2,\,X_3)=-{\alpha_2\alpha_3\over \gamma^2(\gamma+1)}.
\end{eqnarray*}
Since
$$X_3\sim {\cal B}e(\alpha_3,\,\gamma-\alpha_3)$$
we have
$$\b{V}ar(X_3)={\alpha_3(\gamma-\alpha_3)\over \gamma^2(\gamma+1)}.$$
Then
\begin{eqnarray*}
\b{C}ov(Y_1,\,Y_2)&=&{-\alpha_1\alpha_2-\alpha_1\alpha_3
-\alpha_2\alpha_3+\alpha_3(\gamma-\alpha_3)\over \gamma^2(\gamma+1)}\\
&=&{-\alpha_1\alpha_2+\alpha_0\alpha_3\over \gamma^2(\gamma+1)}.
\end{eqnarray*}
It follows that the correlation coefficient is given by
$$\rho(Y_1,\,Y_2)={-\alpha_1\alpha_2+\alpha_0\alpha_3\over
\sqrt{(\alpha_1+\alpha_3)(\alpha_0+\alpha_2)(\alpha_2+\alpha_3)
(\alpha_0+\alpha_1)}},$$
and the correlation is positive if $\alpha_0\alpha_3>\alpha_1\alpha_2$.

\noindent
Suppose now that we want to sample from a bivariate density with beta
marginals, with parameters, respectively, $c_1,\,c_2$ and $c_3,\,c_4$,
and a given positive correlation coefficient $r$. To fit into the previous
framework we set
$$\left\{\begin{array}{lll}c_1&=&\alpha_1+\alpha_3,\\
c_2&=&\alpha_0+\alpha_2,\\
c_3&=&\alpha_2+\alpha_3.\end{array}\right .$$
Then
$$c_4=\alpha_0+\alpha_1=c_1+c_2-c_3,$$
which implies
$$c_1+c_2>c_3.$$
Furthermore
$$r=
{-\alpha_1\alpha_2+\alpha_0\alpha_3\over
\sqrt{c_1c_2c_3(c_1+c_2-c_3)}}.$$
We assume $r>0$.
We solve for $\{\alpha_0,\,\alpha_1,\,\alpha_2,\,\alpha_3\}$, as functions
of $\{c_1,\,c_2,\,c_3,\,r\}$. We get
\begin{equation}
\label{eq:alfa3}
\alpha_3={r\sqrt{c_1c_2c_3(c_1+c_2-c_3)}+c_1c_3\over c_1+c_2}.
\end{equation}
It follows $\alpha_3>0$, as required. And
\begin{equation}
\label{eq:altre}
\left\{\begin{array}{lll}\alpha_1&=&c_1-\alpha_3,\\
\alpha_2&=&c_3-\alpha_3,\\
\alpha_0&=&c_2-c_3+\alpha_3.\end{array}\right .
\end{equation}
All the parameters have to be positive, so we must have
\begin{equation}
\label{eq:restrizioni}
\left\{\begin{array}{lll}c_1-\alpha_3&>&0,\\
c_2-\alpha_3&>&0,\\
c_2-c_3+\alpha_3&>&0.\end{array}\right .
\end{equation}
To determine in the general case the conditions upon which these
restrictions are satisfied is rather cumbersome. We analyze in details
some particular cases.

\smallskip

\noindent
\b{I Case}: $c_1=c_3$. In this case the two marginal have the same
distribution ${\cal B}e(c_1,\,c_2)$.

\noindent
We have
\begin{eqnarray*}
\alpha_3 &=&{r\sqrt{c_1^2c_2(c_1+c_2-c_1)}+c_1^2\over c_1+c_2}\\
&=&{c_1(c_1+rc_2)\over c_1+c_2}.
\end{eqnarray*}
We see that the first two condition in Equation~\ref{eq:restrizioni}
are identical and
\begin{eqnarray*}
c_1-\alpha_3&=&{c_1^2+c_1c_2-rc_1c_2-c_1^2\over c_1+c_2}\\
&=&{c_1c_2(1-r)\over c_1+c_2}\\
&>&0.
\end{eqnarray*}
As for
the third condition we have
\begin{eqnarray*}
c_2-c_3+\alpha_3 &=& {c_2^2+rc_1c_2\over c_1+c_2}\\
&>& 0.
\end{eqnarray*}
We can conclude that in this case we do not have any restrictions on the
parameters.

\smallskip

\noindent
\b{II Case}: $c_2=c_3$.

\noindent
In this case $c_1+c_2-c_3>0$. We have
$$
\alpha_3
={rc_1c_2 +c_1c_2\over c_1+c_2}.$$
\begin{enumerate}
\item
$$c_1-\alpha_3 = c_1(c_1-rc_2).$$
To be greater than zero it requires
$$c_1>rc_2.$$
\item
$$c_3-\alpha_3= c_2-\alpha_3={c_2(c_2-rc_1)\over c_1+c_2}.$$
To be greater than zero it requires
$$c_2>rc_1.$$
The two conditions are then satisfied if
\begin{equation}
rc_2<c_1<{c_2\over r}.
\end{equation}
\item
$$c_2-c_3+\alpha_3=\alpha_3 >0.$$
\end{enumerate}

\smallskip

\noindent
\b{III Case}: $c_1=c_2$.

\noindent
In this case $c_1+c_2-c_3=2c_1-c_3$, so we must have $\delta=2c_1-c_3>0$,
which implies ${c_1\over c_3}>{1\over 2}$.

\noindent
We have
$$
\alpha_3
={rc_1\sqrt{c_3\delta}+c_1c_3\over 2c_1}.$$
\begin{enumerate}
\item
$$c_1-\alpha_3 = {\delta -r\sqrt{c_3\delta}\over 2}.$$
To be greater than zero it requires
$$\delta > r\sqrt{c_3\delta},$$
that is $\delta > r^2c_3$. That means
$$c_1>{c_3(1+r^2)\over 2}.$$
Because
$$c_3\ge {c_3(1+r^2)\over 2},$$
if we have $c_1\ge c_3$ the above condition is always satisfied.
\item
$$c_2-\alpha_3= c_1-\alpha_3,$$
same as before.
\item
$$c_3-\alpha_3= {c_3-r\sqrt{c_3\delta}\over 2}.$$
To be greater than zero we must have
$$c_3>r\sqrt{c_3\delta},$$
that is
$$c_3>r^2\delta,$$
which leads to
$$c_1<{c_3(1+r^2)\over 2r^2}.$$
Putting the two conditions together we have
$${c_3(1+r^2)\over 2}<c_1<{c_3(1+r^2)\over 2r^2},$$
which can be rewritten as
$${1+r^2\over 2}<{c_1\over c_3}<{1+r^2\over 2r^2}.$$
\item
$$c_2-c_3+\alpha_3=c_1-c_3+\alpha_3.$$
We get
\begin{eqnarray*}
c_1-c_3+\alpha_3 &=& {c_1\left [2c_1-c_3+r\sqrt{c_3\delta}\right ]\over 2c_1}\\
&=& {\delta +r\sqrt{c_3\delta}\over 2}\\
&>&0.
\end{eqnarray*}
\end{enumerate}

\smallskip

\noindent
\b{IV Case}: $c_1=c_2=c_3$.

\noindent
Then
$$\alpha_3={c_1(1+r)\over 2}.$$
\begin{enumerate}
\item $$c_1-\alpha_3 ={c_1(1-r)\over 2} >0.$$
\item $$c_3-\alpha_3=c_1-\alpha_3 >0.$$
\item $$c_2-c_3+\alpha_3 = \alpha_3>0.$$
\end{enumerate}

\section{Efficiency with Johnk's method}
If in sampling from the Dirichlet distribution we use a generalization
of Johnk's method (see, for example, \cite[pages 136-137]{chiodi}), that is
\begin{enumerate}
\item generation of 4 independent uniform variates, $U_i,\; i=0,\,\ldots,\,
3$.
\item evaluation of $Z_i=U_i^{1\over\alpha_i},\;
i=0,\,\ldots,\,3$, where $\{\alpha_i\}$ are given by Equation~\ref{eq:alfa3},
Equation~\ref{eq:altre}, under restrictions
as in Equation~\ref{eq:restrizioni}.
\item conditioning upon $S=\sum_{i=0}^nZ_i \le 1$, the variables
$X_i={Z_i\over S},\; i=1,\,2,\,3$ possess the required distribution
${\cal D}(3,\,\alpha_0,\,\alpha_1,\,
\alpha_2,\,\alpha_3)$,
\end{enumerate}
then the efficiency, that is $\prob{S\le1}$, is given by
\begin{equation}
\varepsilon={\prod_{i=0}^3\alpha_i{\m{\Gamma}}(\alpha_i)\over
\left (\sum_{i=0}^3\alpha_i\right ){\m{\Gamma}}
\left (\sum_{i=0}^3\alpha_i\right )}.
\end{equation}
To give an idea we present the following data, referred to the case
$c_1=c_3$. Table 1 refers to the case where $r=0.5$, Table 2 when
$r=0.75$.

\begin{center}
{\em Table 1. Efficiency with Johnk's method \\for selected values of
$c_1=c_3$, $c_2$ and $r=0.50$}

\begin{tabular}{||c|c|c|c|c|c|c|c|c||}\hline\hline
& \multicolumn{8}{c||}{$c_2$}\\ \cline{2-9}
$c_1$&0.23&0.25&0.5&0.75&1.00&2.00&3.00&5.00\\ \hline
0.25&0.907&0.897&0.821&0.763&0.716&0.596&0.526&0.444\\
\hline
0.50&0.833&0.822&0.701&0.612&0.544&0.382&0.298&0.213\\
\hline
0.75&0.778&0.763&0.612&0.506&0.428&0.257&0.179&0.108\\
\hline
1.00&0.734&0.717&0.544&0.428&0.345&0.180&0.113&0.058\\
\hline
1.50&0.668&0.647&0.447&0.323&0.243&0.098&0.050&0.019\\
\hline
2.00&0.619&0.596&0.381&0.257&0.180&0.058&0.024&0.007\\
\hline
3.00&0.552&0.526&0.299&0.179&0.113&0.020&0.007&0.001\\
\hline
5.00&0.472&0.444&0.213&0.108&0.058&0.007&0.001&0.000\\
\hline
\end{tabular}
\bigskip

\bigskip
{\em Table 2. Efficiency with Johnk's method \\for selected values of
$c_1=c_3$, $c_2$ and $r=0.75$}

\begin{tabular}{||c|c|c|c|c|c|c|c|c||}\hline\hline
& \multicolumn{8}{c||}{$c_2$}\\ \cline{2-9}
$c_1$&0.23&0.25&0.5&0.75&1.00&2.00&3.00&5.00\\ \hline
0.25&0.916&0.909&0.843&0.793&0.752&0.646&0.583&0.507\\
\hline
0.50&0.854&0.843&0.734&0.654&0.592&0.441&0.360&0.273\\
\hline
0.75&0.806&0.793&0.654&0.558&0.482&0.315&0.234&0.155\\
\hline
1.00&0.768&0.752&0.592&0.482&0.403&0.232&0.158&0.090\\
\hline
1.50&0.710&0.691&0.503&0.381&0.298&0.138&0.079&0.035\\
\hline
2.00&0.668&0.646&0.441&0.315&0.233&0.088&0.043&0.015\\
\hline
3.00&0.607&0.583&0.361&0.239&0.157&0.043&0.016&0.004\\
\hline
5.00&0.534&0.507&0.273&0.155&0.091&0.015&0.004&0.000\\
\hline
\end{tabular}
\end{center}


\begin{thebibliography}{9}
\bibitem{chiodi} M. Chiodi, \i{Tecniche di simulazione in Statistica},
RCE Edizioni, Napoli, 2000.
\bibitem{kotz} Johnson, N.L. and S. Kotz, \i{Distributions in
Statistics: Continuous Multivariate Distributions}, Wiley, New York, 1972.
\end{thebibliography}
\end{document}